\documentclass[12pt,a4paper]{article}

\usepackage{pdflscape}
\usepackage{amsmath}
\usepackage{amssymb}
\usepackage{latexsym}
\usepackage{srcltx}
\usepackage{graphics}
\usepackage{color}
\usepackage{epsfig}
\usepackage{color}
\usepackage[ruled,vlined]{algorithm2e}

\usepackage{subcaption}

\usepackage{anyfontsize}
\usepackage{tabularx}

\usepackage{graphicx}
\usepackage{multirow}
\usepackage{array}

\usepackage{amssymb, amsmath, amsfonts, mathtools}
\usepackage{esint}
\usepackage{graphicx}
\usepackage{breqn}
\usepackage{float}
\usepackage{caption}
\usepackage{subcaption}
\usepackage{hhline}
\usepackage{multirow}
\usepackage{rotating}
\usepackage{fixltx2e}

\def\la{\langle}
\def\ra{\rangle}

\newcommand{\re}{{\mathbb R}}
\newcommand{\n}{{\mathbb N}}

\newcommand{\cA}{{\mathcal{A}}}

\newcommand{\cE}{{\mathcal{E}}}
\newcommand{\cF}{{\mathcal{F}}}

\newcommand{\cI}{{\mathcal{I}}}
\newcommand{\cJ}{{\mathcal{J}}}

\newcommand{\cV}{{\mathcal{V}}}

\newcommand{\cB}{{\mathcal{B}}}

\newcommand{\cH}{{\mathcal{H}}}
\newcommand{\cS}{{\mathcal{S}}}

\newcommand{\by}{{\boldsymbol y}}
\newcommand{\bx}{{\boldsymbol x}}

\newcommand{\be}{{\boldsymbol e}}

\newcommand{\bv}{{\boldsymbol v}}
\newcommand{\bu}{{\boldsymbol u}}

\newtheorem{theorem}{Theorem}
\newtheorem{prop}{Proposition}
\newtheorem{lemma}{Lemma}
\newtheorem{cor}{Corollary}
\newtheorem{remark}{Remark}

\newtheorem{defi}{Definition}

\date{}

\author{Aleksandar Cvetkovi\' c, \thanks{GSSI (L'Aquila, Italy), 
{e-mail: \tt\small alegzander88@gmail.com}}
Vladimir Yu.~Protasov
\thanks{University of L'Aquila (Italy),  National University High School of Economics, 
Moscow, Russia,  {e-mail: \tt\small
v-protassov@yandex.ru}}}

\title{Maximal acyclic subgraphs and closest stable matrices
\thanks{V.Yu.Protasov acknowledges the support of the Russian
Foundation for Basic Research, projects  no. 19-04-01227 and 20-01-00469. 
  His research was performed within the Basic Research Program
of the National Research University Higher School of
Economics and was supported within the 5-100 Program
for State Support of Leading Universities of the Russian
Federation.
}}

\begin{document}
\maketitle

\begin{abstract}

\bigskip

We develop a matrix approach to the Maximal Acyclic Subgraph 
(MAS) problem by reducing it to finding the closest nilpotent matrix to the 
matrix of the graph. Using recent results on the closest Schur stable systems  
and on minimising the spectral radius over special sets of non-negative matrices 
we obtain an algorithm for finding an approximate  solution of MAS. 
Numerical results for graphs from $50$ to $1500$ vertices 
demonstrate its fast convergence and give the rate of approximation 
in most cases 
larger than $0.6$.  The same method gives the precise solution 
for the following   weakened version of 
 MAS:  find the minimal $r$ such that the graph can be made acyclic 
by cutting at most~$r$ incoming edges from each vertex. 
Several modifications,  
when each vertex is assigned with its own maximal number~$r_i$ of cut edges,  when some of edges are ``untouchable'',  are also considered. Some applications are discussed.

\bigskip

\smallskip

\begin{flushright}
\noindent  \textbf{AMS 2010} {\em subject
classification: 05C50, 15B48, 90C26}
\end{flushright}

\end{abstract}
\vspace{1cm}

\begin{center}
\large{\textbf{1. Introduction}}
\end{center}
\vspace{1cm}

Let $G = (\cV, \cE)$ be a directed graph with the set of vertices 
$\cV = \{g_1, \ldots , g_n\}$ and with the set of edges~$\cE$. 
The vertices   will be identified with 
the corresponding numbers~$\{1, \ldots , n\}$.
 The {\em Maximal Acyclic Subgraph} (MAS) of $G$ is a graph  
 $\widehat G = (\cV, \widehat \cE)$ such that $\widehat \cE \subset \cE$, 
 $\widehat G$ has no cycles,  
and the number of edges $|\widehat \cE|$ is maximal.  The MAS problem (in short MAS) is to find a maximal acyclic 
subgraph of a given graph~$G$. This problem  was included 
by R.Karp in 1972 in his list of 21 NP-complete problems~\cite{K}.  

Most of methods  in the literature usually find approximate solutions to the MAS problem.
We say that a  subgraph $G'$ of $G$ gives an approximate solution with an {\em approximation factor}~$\delta \in  (0, 1)\, $ if $G$ is acyclic and has at least $\, \delta \, |\widehat \cE|$ edges. 
It is not known if efficient algorithms exist to obtain an approximate solution
with a factor $\delta > \frac12$, 
and most likely the answer is negative~\cite{GMR}.

It is interesting that  the simplest algorithm for an approximate solution of MAS gives the best known result $\delta = \frac12$. 
Indeed, take an arbitrary renumbering of vertices of $G$ 
(a permutation~$\sigma$ of the set $\{1, \ldots , n\}$), then the
set of forward edges in the new ordering and the set of backward edges are both acyclic and 
 one of them contains at least $\frac12 \, |\cE| \, \ge \, \frac12 \, |\widehat \cE|$ edges. Hence, every permutation gives 
an approximate solution to 
MAS with the factor  $\delta \ge \frac12$. 
In fact, the exact solution of MAS is equivalent to finding an optimal  permutation~$\sigma$ with the largest number of forward edges.  
See~\cite{BS, CMM, GMR, HR} for  algorithms of approximate solutions and for discussions of the complexity.

We develop a linear algebraic approach to MAS. First we  reformulate this problem as 
finding the closest (in the Frobenius norm) non-negative matrix~$X$  to a given  matrix~$A$ with the condition that $\rho(X) = 0$ (the spectral radius  is equal to zero). 
One can note that this problem is related to  finding the 
closest Schur stable matrix to a given matrix~$A$. That stabilisation problem
is formulated in the same way but with the condition $\rho(X) = 1$. 
Some methods for the Schur stabilization were elaborated in recent literature, 
but all of them find only local minima, which does not give much for the MAS 
problem. Nevertheless, replacing the Frobenius norm by the $L_{\infty}$-norm, 
we obtain the Schur stabilization problem that can be effectively solved~\cite{NP2}.   
In graph terms this replacement leads to the following problem: {\em find the minimal number $r$ such that one can make the graph acyclic 
by cutting at most $r$ incoming edges from each vertex.} In Section~3 we present an 
algorithm that  solves this problem. For graphs with 100 vertices, finding the solution   takes  usually a few seconds in a standard laptop (detained characteristics are given in Section~6), for 1000 vertices it takes about
10  minutes.  See Section~6 for a discussion of numerical aspects. The algorithm is based on minimising  the Perron eigenvalue 
on special sets of matrices. In Section~4 we formulate several generalisations of our  problem which can be solved by  similar methods. In Section~5 we come back to 
the (classical) MAS problem. We show that a slight 
modification of our method can be used to finding an approximate solution of MAS. 
The  rate of approximation is estimated from below 
by the value $\, \gamma \, = \, \frac{|\cE'|}{|\, \cE \, |}$, where~$\cE'$
is the set of edges of the obtained acyclic subgraph. 
Since  $\delta \, = \, \frac{|\cE'|}{|\, \widehat \cE \, |}  \, \ge \, 
 \frac{|\cE'|}{|\, \cE\, |}$, 
we see that~$\delta \, \ge \, \gamma$. 
In the numerical experiments, our method  gives the rate of approximation 
better than $\gamma \approx 0.6$. 
In Section~6 we present numerical results and discuss the complexity issue, in Section~7 
we consider some applications. 
 \smallskip 
 
Throughout the paper, by the 
graph of a non-negative $n\times n$ matrix~$A$ we mean the graph 
$G$ with $n$ vertices $\{1, \ldots , n\}$ such that there is an edge from a
vertex $i$ to a vertex $j$ if and only if $A_{ji}>0$.
We use a standard componentwise ordering of real matrices: 
$A \ge (>) B$ if $A_{ij} \ge  (>) B_{ij}$ for all $i,j$; the same for vectors;  
$\rho(A)$ denotes the spectral radius of a matrix~$A$, i.e., the largest modulus of its eigenvalues. By the Perron-Frobenius 
theorem~\cite[chapter 8]{G}, if $A \ge 0$, then there is an  eigenvalue 
of $A$ denoted 
$\lambda_{\max}$ equal to the spectral radius of~$A$ and there is  a non-negative eigenvector of~$A$ 
corresponding to~$\lambda_{\max}$. This eigenvalue
and the corresponding eigenvector are referred to as {\em leading}. 
The vectors will be denoted by bold letters and their components by standard 
letters: $\bx = (x_1, \ldots , x_n)$.  
Matrices (vectors) with entries $0$ and $1$
are Boolean. In particular, an adjacency matrix of a graph 
($A_{ji} = 1 \, \Leftrightarrow \, $ there is an edge from $i$ to $j$; 
otherwise $A_{ji} = 0$) is Boolean.

We use the notation $\cV = \{1, \ldots , n\}$. This set will be identified with 
the set of vertices of the graph~$G$. For an arbitrary finite set $K$, 
$|K|$ denotes its cardinality. 
Let $X$ be an $n\times n$ matrix. For a subset $\cS \subset \cV$, we denote by 
$X|_{\cS}$ the $|\cS|\times |\cS|$
principal submatrix of $X$, which is a restriction of $X$ to the index 
set $\cS$. By $\bx|_{\cS}$ we denote the restriction of the vector $\bx \in \re^n$
to the $|\cS|$-dimensional space of vectors with all non-zero components in~$\cS$.

The Frobenius norm is the Euclidean norm on the set of matrices understood as vectors in the $n^2$ dimensional space: $\|X\| = \bigl[ \sum_{i,j} |X_{ij}|^2\bigr]^{1/2}$.
We always use this matrix norm if another norm is not specified.

\vspace{1cm}

\begin{center}
\large{\textbf{2. The MAS problem vs the closest stable matrix problem}}
\end{center}
\vspace{1cm}

The MAS problem  can be formulated with spectra of graphs 
by applying  the well-known fact: {\em a directed graph is acyclic if and only if 
the spectral radius of its adjacency matrix is zero}. Indeed, if $\rho(A) = 0$, 
then $A^n = 0$ (this follows easily if we write~$A$ in its Jordan normal form).  
Hence the graph~$G$ has no walks of length~$n$. Therefore it cannot have cycles, otherwise a cycle generates walks of all lengths. Conversely, if $G$ is acyclic, then it cannot 
have long walks, and hence some high power of~$A$ is equal to zero, which implies~$\rho(A) = 0$.   

The following simple observation presents 
the MAS problem in the linear algebraic terms. 
\begin{prop}\label{p.10}
The MAS problem is equivalent to the following problem: 
given a graph~$G$ with an adjacency matrix~$A$, find the closest 
in Frobenius norm non-negative matrix~$X$ with zero spectral radius: 
\begin{equation}\label{eq.prob-frob}
\left\{
\begin{array}{l}
\|X-A\| \ \to \ \min\,  \\
X \ge 0\, , \ \rho(X) = 0 \, .
\end{array}
\right.
\end{equation}
Every solution $X$ of the problem~(\ref{eq.prob-frob}) is a Boolean matrix whose graph solves the MAS problem
for~$G$. 
\end{prop}
{\tt Proof.} The matrix $A$ is Boolean. Hence, if we replace all
strictly  positive entries of $X$ by ones, then the spectral radius remains zero but the distance $\|X-A\|$
becomes smaller. If for each pair $(i,j)$, we set $X_{ij} = 0$ whenever 
$A_{ij} = 0$, we obtain the same effect. Consequently, problem~(\ref{eq.prob-frob}) 
can be reduced to the set of Boolean matrices $X$ such that $X\le A$ i.e., to matrices of subgraphs of $G$. In this case, $\|X-A\|$ is equal to the square root 
of the number of edges cut to obtain $X$ from $A$. On the other hand, $\rho(X)= 0$ 
precisely when the graph of~$X$ is acyclic. Thus, 
problem~(\ref{eq.prob-frob}) is equivalent to finding the acyclic subgraph 
obtained from $A$ by cutting the smallest number of edges, i.e., to 
the MAS problem.

{\hfill $\Box$}
\smallskip

Problem~(\ref{eq.prob-frob}) looks similar to the {\em closest stable matrix} problem. We formulate only the non-negative version of this problem. 
Given a matrix $A \ge 0$,  find the closest matrix such that~$\rho(X) = 1$:  
\begin{equation}\label{eq.prob-stab}
\left\{
\begin{array}{l}
\|X-A\| \ \to \ \min\,  \\
X \ge 0\, , \ \rho(X) = 1 \, . 
\end{array}
\right.
\end{equation}
This is the same problem~(\ref{eq.prob-frob})
but with $\rho = 0$ replaced by $\rho=1$. The closest stable matrix problem 
have been studied in many recent works due to its application in 
dynamical and controlled systems, mathematical economics, population dynamics, etc. 
(see~\cite{By, GM, MMS, ONVD} for the problem with general matrices and 
\cite{An, NP2, GP2} for non-negative matrices). A question arises if  the 
methods for the closest stable matrix problem can be useful for MAS? 
The main difficulty is that all known methods find only locally stable matrices. 
Moreover, problem~(\ref{eq.prob-frob}) may have an exponential number of local minima~\cite{GP2}. It turns out, however, that in another matrix norm, the closest stable matrix can be efficiently found. For example, in  
the $L_{\infty}$-norm: $\, \|X\|_{\infty}  \, = \, \max_{i=1, \ldots , n}\sum_{j=1}^n |X_{ij}|$. An algorithm presented in~\cite{NP2} finds the global minimum 
in the problem~(\ref{eq.prob-stab}) with $L_{\infty}$-norm quite fast even in high 
dimensions (several thousands). Reformulating the problem in a new norm in term of graphs we come to the following  variant of MAS:  
\smallskip 

\noindent \textbf{The max-MAS Problem.} {\em Find the minimal integer $r$
such that a given graph $G$ can be made acyclic by cutting at most 
$r$ incoming edges from each vertex}. 
\smallskip 

Suppose a graph $G$ is made acyclic by cutting some edges. 
If $c_i$ denotes 
the number  of cut incoming edges from the $i$th vertex, then   
 the MAS problem minimises the sum $\sum_{i=1}^n c_i$ while the max-MAS 
problem minimises $\max_{i=1, \ldots , n}\, c_i$. This justifies our terminology. 
In the next section we will see that 
this problem can be efficiently solved by the method of
finding the closest stable matrix~\cite{NP2}. Note that in spite of similarity of problems
(\ref{eq.prob-frob}) and  (\ref{eq.prob-stab}) they are different, 
since the set of matrices with $\rho(X) = 1$ and with $\rho(X) = 0$
have different structures. Nevertheless, the max-MAS problem can be solved 
by the method based on the algorithms for solving the problem~(\ref{eq.prob-stab}) in $L_{\infty}$-norm.
Moreover, the same method can solve several generalisations of max-MAS. For example,  
when it is allowed to cut at most~$r_i$ incoming edges from the vertex~$i$, where 
$r_1, \ldots , r_n$ are given numbers. Some of those numbers may be zeros, in which case 
the corresponding vertex is ``untouchable''. This and other generalizations are solved in Section~4. Then, in Section~5 we suggest an approximate solution of the 
classical MAS problem based on the presented algorithms. Numerical results and the 
complexity issue are discussed in Section~6. In Sections 7 we discuss possible applications  
and consider the application to small-world networks in detail. 

\newpage

\begin{center}
\large{\textbf{3. Algorithmic solution for the max-MAS problem }}
\end{center}
\vspace{1cm}

We first reformulate the MAS and max-MAS problems in terms of optimising 
the spectral radius of a non-negative matrix. Then we establish a relation 
between those problems and the problem of stabilization of a positive 
linear system. Then we apply the technique of stabilisation 
to construction of acyclic graphs. We will see that the max-MAS problem can 
be solved completely by this approach. Approximate solutions for MAS will be considered in Section~5. 
\vspace{1cm}

\begin{center}
\textbf{3.1. Spectral formulation of the max-MAS problem}
\end{center}
\vspace{1cm}

The max-MAS problem is  formulated in linear algebraic terms as follows: 
\begin{equation}\label{eq.max-mas}
\left\{
\begin{array}{l}
\|X-A\|_{\infty} \ \to \ \min\,  \\
X \ge 0\, , \ \rho(X)  = 0\, , 
\end{array}
\right.
\end{equation}
where $A$ is a given Boolean $n\times n$ matrix. 
The equivalence of this problem to max-MAS is proved in the same way as 
in Proposition~\ref{p.10}. This problem, in turn, is reduced to 
\begin{equation}\label{eq.prob10}
\left\{
\begin{array}{l}
\rho(X) \ \to \ \min\,  \\
X \ge 0\, , \ \|X-A\|_{\infty} \ \le \ r \, . 
\end{array}
\right.
\end{equation}
Indeed, if we are able to solve~(\ref{eq.prob10}) for every 
integer $r$, then the minimal possible $r$ for which 
there is a non-negative matrix $X$ such that $\|X-A\|_{\infty} \le r$ and $\rho(X) = 0$
is found merely by the integer bisection in~$r$. 

It is more convenient to consider problem~(\ref{eq.prob10}) in slightly 
different terms. For the $i$th row of the matrix~$A$, 
we denote $\cB(A_i, r) \, = \, \bigl\{ \bx \in \re^n_+ \ \bigl| \ 
 \|A_i - \bx\|_1 \le r\}$, where $\|\by\|_1 = \sum_{j=1}^n |y_j|$ is the $L_1$-norm. So, $\cB(A_i, r)$ is the $L_1$-ball of radius $r$  centered at $A_i$
 and intersected with the positive orthant. 
 Then $\cB(A, r)$ denotes the set of matrices with the $i$th row from $\cB(A_i, r), \, i = 1, \ldots , n$. Clearly, $\cB(A, r)$ is a convex polyhedron in the set of matrices. 
 
 Now we focus on the following problem equivalent to~(\ref{eq.prob10}): 
 for a given Boolean matrix~$A$ and for $r\in \n$, solve 
\begin{equation}\label{eq.prob20}
\left\{
\begin{array}{l}
\rho(X) \ \to \ \min\,  \\
X \, \in \, \cB(A, r)  \, . 
\end{array}
\right.
\end{equation}
The solution $X$ is always a Boolean matrix. If we find  the smallest non-negative 
integer~$r$ for which this problem has a solution with the objective function value equal to zero, then the graph of this  
solution~$X$ 
solves the max-MAS problem. 
\vspace{1cm}

\begin{center}
\textbf{3.2. Minimisation of spectral radius over product families of non-negative matrices }
\end{center}
\vspace{1cm}

The problem~(\ref{eq.prob20}) can be efficiently solved for each~$r$. 
This is because $\cB(A, r)$ is a {\em product set} of matrices.
Optimising the spectral radius over product sets has been 
investigated in various contexts~\cite{NP1, P16, NP2, 
Akian1}. Let us have arbitrary compact sets $\cF_i \in \re^n_+$. 
The set of matrices $\cF = \{X \ | \  X_i \in \cF_i, \, i = 1, \ldots , n\}$
is called a product set. Every matrix from a product set 
is composed with rows chosen from the sets $\cF_i$ arbitrarily and 
independently. The sets $\cF_i$ are usually referred to as {\em uncertainty sets}. 
All methods of minimising the spectral radius over product families
are actually based on the well-known formula 
\begin{equation}\label{eq.cw}
\rho(X)\ =\ \sup\, \bigl\{ \, \lambda \  \bigl| \
\exists \, \bu \ge 0 \, , \, \bu \ne 0\ : \  X\, \bu \ \ge \ \lambda \, \bu\, \bigr\}\, ,  
\end{equation}
which holds  for every non-negative matrix~\cite{BP}. 
If we minimise the spectral radius of a matrix over a product set,  
this formula allows us to treat all rows separately: 
\begin{equation}\label{eq.cw-p}
\min_{X \in \cF}\, \rho(X)\ =\ \sup\, \bigl\{ \, \lambda \ \bigl| \
\exists \, \bu \ge 0 \, , \ \bu \ne 0\ : \  \la X_i \, , \,  \bu \ra \ \ge \ \lambda \, u_i\, 
,  X_i \in \cF_i, \, i = 1, \ldots , n\bigr\}\, . 
\end{equation}
If all the uncertainty sets $\cF_i$ are finite, then there are methods of applying linear programming to this problem. 
Their complexity, however, depends on the cardinalities  of the uncertainty 
sets~\cite{PY}.  
In our case the sets 
$\cF_i = \cB(A_i, r)$ are polytopes with exponentially many vertices. 
If the row~$A_i$ contains $m_i$ ones and $m_i > r$, 
then the number of vertices is ${m_i \choose r}$. Therefore, the linear
problem has an exponential number of inequalities, which makes this approach 
inefficient 
for the problem~(\ref{eq.prob20}). 

That is why we make use of another approach based on a 
recursive relaxation   scheme suggested in~\cite{Akian1, NP2}
and called {\em the greedy algorithm}.  
 A crucial point of this approach is the following fact: 
\smallskip

\noindent \textbf{Theorem A.}~\cite{NP1} {\em 
A matrix $X \in \cF$ has the minimal spectral radius over a product set~$\cF$
if and only if $X$ possesses a leading eigenvector~$\bv$ such that 
$\la X_i , \bv\ra \, = \, \min\limits_{\bx \in \cF_i}\la \bx , \bv\ra$
for each $i = 1, \ldots , n$. }
\smallskip 

 The relaxation scheme works as follows. If in $k$th iteration we have a matrix $X^{(k)} \in \cF$, 
then we compute its leading eigenvector $\bv^{(k)}$ and then,  for every $i=1, \ldots , n$, we replace the $i$th row of $X^{(k)}$ by the element  $X_i \in \cF_i$
which minimises the scalar product $\la X_i, \bv^{(k)} \ra$. We obtain 
$X^{(k+1)}$, etc. Thus, in each iteration we replace every row of the matrix by 
the optimal row in the corresponding uncertainty set, i.e., by the 
row making the shortest projection with the leading eigenvector. 
Formula~(\ref{eq.cw}) implies that this scheme is a relaxation: 
 $\rho(X^{(k)})$ is decreasing in 
$k$, maybe non-strictly. Under some ``positivity-like'' assumptions on the 
uncertainty sets $\cF^{(k)}$, 
the spectral radius strictly decreases, and hence the solution is found 
in finite time. This is true, for instance, when all vectors from all the sets $\cF^{(k)}$
are strictly positive, or when each matrix from the product family $\cF$
is irreducible (does not
have a nontrivial invariant coordinate subspace, i.e., a
subspace spanned by several vectors of the canonical
basis). 

This scheme, however, has a serious disadvantage: without those 
  ``positivity-like'' assumptions the algorithm may cycle. Moreover, for sparse matrices, the cyclicity happens quite often.  A modified 
  greedy algorithm which never cycles was presented in~\cite{NP2}. According to 
  numerical experiments, that algorithm has a very fast convergence and finds the 
  munimum within a few iterations even in high dimensions~\cite{CP}. Of course, if some of the sets 
  $\cF_i$ are of large cardinality, one iteration may take long, since it requires 
  computing of $|\cF_i|$ scalar products. This occurs, in particular, 
  in our problem~(\ref{eq.prob20}), where each set $\cF_i = \cB(A_i, r)$
  has an exponential number of extreme points.  In this section we overcome this difficulty 
  and modify the greedy algorithm specially for these sets.  To this end we begin with introducing some further notation.    
\smallskip 

An arbitrary non-negative vector  $\bv \in \re^n_+$ 
defines an ordering on the set $\cV = \{1, \ldots , n\}$ by the values 
of components of~$\bv$: $\ i \, \ge_{\bv} \, j\, $ if 
$v_i \ge v_j$. For every $r\le n$, we can consider the set of $r$ largest elements of 
$\cV$ in this ordering. This set may not be unique if~$\bv$ have some equal components. In this case we take any 
set of $r$ largest elements. Similarly, for any subset $\cS \subset \cV$, 
we consider the set~$\cS'$ of $r$ largest elements of~$\cS$. We say that the set~$\cS'$    
{\em occupies the $r$ largest components of $\bv$}. If $|\cS| < r$, then 
we say that the whole set $\cS' = \cS$ occupies the~$r$ largest components of $\bv$
(although it contains fewer number of elements).

\begin{defi}\label{d.20}
Let $A$ be a Boolean  $n\times n$ matrix and $r\in \n$. 
Let also  $X$ be a  Boolean $n\times n$ matrix
and $\bv \in \re^n_+$ be a vector. 
A row $X_i$ is minimal 
in the ball $\cB(A_i,r)$ with respect to~$\bv$
if the set of zeros of $X_i$ on the set $\, {\rm supp}\, A_i$ 
occupies $r$ largest components of $\bv$ on that set.  
\end{defi}
Thus, if $\|A_i\|_1 \, > \, r$, then the minimality 
of $X_i$ means that
in the index set $\, {\rm supp}\, A_i$, there are 
$r$ indices corresponding to the maximal (in this set)  components of the vector 
$\bv$ on which all components of $X_i$ are zeros. 
If $\|A\|_1 < r$, then the  minimality means simply that 
 $X_i = 0$ on the whole set $\, {\rm supp}\, A_i$.  

It is easy to find the minimal row in the ball~$\cB(A_i,r)$  with respect to 
a given vector~$\bv$. One needs to order the index set $\, {\rm supp}\, A_i$ by the values of the components 
of~$\bv$ and take the $r$ largest (in this order) indices. 
Denote this set by $\cJ$. If $|{\rm supp}\, A_i| \le r$, then put 
$\cJ = {\rm supp}\, A_i$. Then the minimal row~$X_i$ is defined as follows: 
$X_{ij} = 0$ if $j\in \cJ$ and $X_{ij} = A_{ij}$ otherwise.

The minimal row possesses the shortest possible projection 
to the vector~$\bv$ among all elements of~$\cB(A_i, r)$ as the following lemma 
asserts.   
\begin{lemma}\label{l.10}
A row $X_i$ is minimal with respect to a vector~$\bv$
in the ball $\cB(A_i,r)$  precisely when 
\begin{equation}\label{eq.min1}
\la X_i, \bv \ra \  = \ \min_{\bx \in \cB(A_i,r)}\la X_i, \bv \ra \, . 
\end{equation}
\end{lemma}
{\tt Proof}. Assume $|{\rm supp}\, A_i| \, > \, r$. 
If the set of zeros of $X_i$ does not occupy 
the $r$ largest components of $\bv$ on ${\rm supp}\, A_i$, 
then there are numbers $j,k \in {\rm supp}\, A_i$ such that 
$v_{j} > v_k$ and $v_j = 1, v_k = 0$. 
Interchanging those components we reduce the scalar  
product $\la X_i, \bv \ra$, which is a contradiction. 
Conversely, assume the set of zeros of $X_i$ occupies  
$r$ largest components of $\bv$ on ${\rm supp}\, A_i$, 
but the minimal scalar product is smaller than for $X_i$ and is attained 
at some $\bx \in \cB(A_i,r)$. It can be assumed that 
 $\bx$ is an extreme point of the ball 
 $\cB(A_i,r)$, i.e. a Boolean vector. Hence 
 $\la A_i - \bx, \bv \ra$ does not exceed the sum of $r$ largest on the 
 set ${\rm supp}\, A_i$ components of~$\bv$, i.e., 
 does not exceed $\la A_i - X_i, \bv \ra$. Hence 
 $\la \bx, \bv \ra \, \ge \, \la X_i, \bv \ra$.
If $|{\rm supp}\, A_i| \, \le \, r$, then 
the minimal (coordinatewise) element of the ball $\cB(A_i, r)$ is the origin, 
for which the minimal scalar product is attained.   

{\hfill $\Box$}
\smallskip

\begin{remark}{\em 
If $\bv$ possesses some zero components, then 
the definition of the minimal row can be reduced to the 
support of~$\bv$. We reduce all vectors to the set $\cS = {\rm supp}\, \bv$
and do not pay attention to other components. Then all the minimal vectors stay minimal after this reduction. 
}
\end{remark}

\begin{theorem}\label{th.10}
Let $A$ be a Boolean matrix and $r\in \n$. Let 
the problem~(\ref{eq.prob20})
reach its global minimum at some Boolean matrix $X$. 
Then this matrix~$X$ is characterised by the property: 
there exists a leading eigenvector~$\bv$ of $X$
such that each row $X_i$ is minimal in the ball $\cB(A_i, r)$ with respect to~$\bv$.  
\end{theorem}
{\tt Proof} follows by combining Lemma~\ref{l.10} and Theorem~A for the sets 
$\cF_i = \cB(A_i, r)$ and for the leading eigenvector $\bv$ of $X$. 

{\hfill $\Box$}
\smallskip

If the eigenvector~$\bv$ is sparse, then it makes sense to 
reduce all vectors to the set~$\cS = {\rm supp} \, \bv$ in the spirit of 
Remark~\ref{r.30}. 
Indeed, the scalar  product $\la X_i , \bv \ra$ depends only 
on entries of the vector $X_i$ on the set~$\cS$. Hence, if $\bv$ is not strictly positive, the criterion of 
Theorem~\ref{th.10} can be reduced to the set $\cS$:  
\begin{cor}\label{c.10}
The criterion of Theorem~\ref{th.10} characterising absolute minima in the problem~(\ref{eq.prob20})
 can be written in the following form. 
Let $\cS \, = \, {\rm supp} \, \bv$. Denote by $\bv', X_i'$, and $A_i'$
the restrictions of those vectors to the set $\cS$. 
Then the matrix $X$ is a solution of the problem~(\ref{eq.prob20})
if and only if 
there exists a leading eigenvector~$\bv$ of $X$
such that for every $i \in \cS$ each row $X_i'$ is minimal  in the ball 
$\cB(A_i', r)$ with respect to $\bv'$.  
\end{cor}  
If $\bv > 0$, then Corollary~\ref{c.10}
coincides with Theorem~\ref{c.10}. If $\bv$ has some zero components, then 
the criterion of Corollary~\ref{c.10}  is simpler in practice since 
it involves only the sumbatrix~$X' = X|_{\cS}$. 

\vspace{1cm}

\begin{center}
\textbf{3.3. Some auxiliary facts on non-negative matrices}
\end{center}
\vspace{1cm}

A non-negative matrix $A$ is called {\em irreducible} if it does not
have a nontrivial invariant coordinate subspace, i.e., a
subspace spanned by some elements~$e_i$ of the canonical
basis. A matrix is irreducible if and only if its graph 
is strongly connected. 
Reducibility  means that there is a proper nonempty subset 
$\Lambda \subset \cV$ such that for each $i \in \Lambda$,
the support of the $i$th column of $A$ is contained in $\Lambda$.
It is well-known (e.g.~\cite{BP}) that 
an irreducible matrix has a simple leading eigenvalue.
The converse is not true: a matrix with a simple leading
eigenvalue can be reducible.

For every matrix $A\ge 0$,  there exists  a suitable
permutation~$P$ of the basis of $\re^n$, after which $A$
gets  a block upper triangular form with $q \ge 1$
diagonal blocks $A_j$ called the {\em Frobenius factorization}:
\begin{equation}\label{eq.frob}
P^{-1}AP \quad = \quad \left(
\begin{array}{cccccc}
A_1 & * &  \ldots &  * \\
0 & A_2&  * & \vdots \\
\vdots & {} &  \ddots  & * \\
0 &  \ldots &  0 & A_q
\end{array}
\right)\ .
\end{equation}
For each $j = 1, \ldots , q$, the matrix $A_j$ in the
$j$th diagonal block is irreducible. Any non-negative
matrix possesses a unique Frobenius factorization up to a
permutation of blocks (see \cite[chapter 8]{G}).

Let $A$ be an $n\times n$ non-negative matrix. Its leading
eigenvector $\bv$ is called {\em minimal} if there is no
other leading eigenvector that possesses a strictly
smaller (by inclusion) support. 
A minimal eigenvector is unique up to normalisation if and only if the 
geometrical multiplicity of the leading eigenvector is one. 
A minimal leading
eigenvector can be found by Frobenius
factorization~(\ref{eq.frob}). It suffices to
take the smallest $m$ such that $\rho(A_m) = \rho(A)$ i.e., the ``highest''
diagonal block with the maximal spectral radius.  Then consider the 
leading eigenvector~$\bu$ of the principal submatrix of $A$ 
that consists of blocks  $A_1, \ldots , A_m$. Complement 
this vector by zeros to an $n$-dimensional vector and obtain a  
minimal eigenvector~$\bar \bu$ 
of the matrix $P^{-1}AP$. Respectively, the vector $P\bar\bu$ is the minimal leading eigenvector of $A$. 

Can the minimal leading eigenvector be strictly positive, i.e., possess a full support? Of course. This case  is characterized by the
following property: 
\begin{prop}\label{p.30}\cite{NP2}
If  a matrix $A\ge 0$ has a strictly
positive minimal leading eigenvector~$\bv$, then 
its leading eigenvalue is simple. In the Frobenius factorisation~(\ref{eq.frob}) the spectral radius 
of $A_q$ is equal to $\rho(A)$ and the spectral radii of all the 
other blocks are smaller than~$\rho(A)$. 
\end{prop}
In case of positive minimal leading eigenvector we can define 
the notion of a {\em basic set}, which is needed in our algorithm. 
\begin{defi}\label{d.30}
Suppose a matrix~$A$ has a strictly
positive minimal leading eigenvector. 
Then the basic set of $A$ is the support of the leading eigenvector 
of the matrix~$A^T$.  
\end{defi}
If the minimal leading eigenvector is strictly positive, then 
by Proposition~\ref{p.30}, the leading eigenvalue is simple. 
Then so is the leading eigenvalue of of the transposed matrix.  This implies the 
correctness  of Definition~\ref{d.30}. Another consequence of 
Proposition~\ref{p.30} is that the basic set can be found without computation 
of the leading eigenvector of~$A^T$, provided the Frobenius factorisation
is available. 
\begin{prop}\label{p.50}\cite{NP2}
Suppose a matrix~$A$ has a strictly
positive minimal leading eigenvector;  
then the basic set of $A$ is the set of indices which 
after the permutation~$P$ correspond to the the last block~$A_q$ 
in the Frobenius factorisation~(\ref{eq.frob}). 
\end{prop}
Thus, if~$A$ has a strictly
positive minimal leading eigenvector, then the basic set  of~$A$ 
can be found as follows. We consider the permutation~$P$
of the basis vectors which takes $A$ to its  Frobenius form~(\ref{eq.frob}). 
Then we take the last block~$A_q$ in this factorisation and denote by $d$ its size. 
Then the basic set is $\, \bigl\{ k \ \bigl| \ 
e_{k} \, = \,  P^{-1}e_{j}, \, j = n-d+1, \ldots , n \bigr\}$, 
where $\{e_j\}_{j=1}^n$ is the canonical basis in~$\re^n$. 
In particular, if $A$ has a Frobenius form, i.e., if $P$ is the identity permutation, 
then the basic set is $\{n-d+1, \ldots , n\}$. 
\smallskip 

{\tt Proof of Proposition~\ref{p.50}.} By Proposition~\ref{p.30} the matrix $A$ and hence 
$P^{-1}AP$ has a simple leading eigenvalue which is located in the 
last block $A_q$ of factorisation~(\ref{eq.frob}). Hence the leading eigenvector of 
$[P^{-1}AP]^T$ has its support in the last block. On the other hand, 
$[P^{-1}AP]^T \, = \, P^{-1}A^TP$, from which the proposition follows. 

{\hfill $\Box$}
\smallskip

Now the preliminary work is done and we are ready to present the main result. 

\vspace{1cm}

\begin{center}
\textbf{3.4. The algorithm for  minimising \\
the spectral radius over an $L_1$-ball}
\end{center}
\vspace{1cm}

We describe the algorithm for solving problem~(\ref{eq.prob20})
of minimising 
the spectral radius over an $L_1$-ball~$\cB(A, r)$, by which we 
find the closest stable matrix on a product family.  The radius $r$ is supposed to be 
a natural number. 
\bigskip

\noindent \textbf{Notation for Algorithm~1}. 
Let $\cV = \{1,\ldots, n\}$. For a matrix $X$, vector $\bx$, and a subset $K \subset \cV$, we denote by $X|_{K}$ and $\bx|_{K}$ the restrictions of the matrix and of the vector to the index set $K$ (see Introduction). We use the 
notions of minimal leading eigenvector (Definition~\ref{d.20}) and of the 
basic set (Definition~\ref{d.30}). We use an obvious fact that if $\bv$ is a minimal leading eigenvector of a matrix~$X$ and $\cS \, = \, {\rm supp}\, \bv$, then 
  the matrix $X|_{\cS}$ has a strictly positive minimal leading 
  eigenvector $\bv|_{\cS}$. 

\vspace{7cm}

\begin{algorithm}
\KwData{A Boolean $n\times n$  matrix~$A$, a number $r \in \n$}
\KwResult{A matrix~$\widehat X$
which is a solution of problem~(\ref{eq.prob20}).   }
\Begin{
 Set  $X^{(1)} = A$, $k = 1$. 
 \smallskip 
 
\nl \textbf{(*)} {\tt $k$th iteration.} We have a Boolean $n\times
n$ matrix $X^{(k)}$. 

\nl \eIf{$\rho(X^{(k)}) < 1$}{\textbf{STOP}. 
Algorithm~1 terminates. Denote $\widehat X = X^{(k)}$ and \textbf{Return}\; }{
\nl 
Denote $X = X^{(k)}$.  Compute a minimal leading eigenvector
$\bv$ of $X$ (take any of them, if there are several ones),  set
$\cS = \cS^{(k)} = {\rm supp} \, \bv$, $X'=X|_{\cS}, \, \bv' = \bv|_{\cS}$\;
\smallskip 

\nl \textbf{(**)} {\tt Main loop.} 
{\em Input}: a triple $(X, \bv, \cS)$, 
where $X$ is Boolean $n\times n$ matrix, $\bv$ is a minimal 
leading eigenvector of~$X$ and $\cS  = {\rm supp}\, \bv $. 
{\em Output}: a  new triple $(X, \bv, \cS)$.

Let  $\cH \subset \cS$ be 
the basic set of~$X'$ and $\cI$ be the set of indices $i\in \cS$ 
for which the row $X_i'$ is minimal in the ball 
$\cB(A_i', \bv')$ with respect to the vector $\bv'$.

\nl \eIf{$\cI = \cS$}{\textbf{STOP}. 
Algorithm~1 terminates. Define the $n\times n$ matrix $\widehat X$ as follows: 
 ${\widehat X}_i= X_i$ for $i \in \cS$ and 
${\widehat X}_i = A_i$ for $i \notin \cS$. \textbf{Return}\; }{
\nl Define the next matrix~$\tilde X$ as follows:

If $i \in \cI$ or $i \notin \cS$, then  $\tilde X_i = X_i$; 
 
Otherwise,  if $i \in \cS\setminus \cI$, then 
$\tilde X_i|_{\cS}$ is the minimal 
row in the ball $\cB(A_i', r)$ with respect to $\bv'$ and 
$\tilde X_{ij} = A_{ij}$ for all $j\notin \cS$. 

Set  $\tilde X' = \tilde X|_{\cS}$  \;  
\nl \eIf{$\cH \subset \cI$}{
  $\rho(\tilde X') = \rho(X')$, 
the leading eigenvalue of $\tilde X'$ is simple 
\; 
\nl We compute the minimal leading eigenvector
$\tilde \bv'$ of $\tilde X'$. The set $\cS$ is not changed. 
Set $X = \tilde X, \bv = \tilde \bv\,$ \;
\nl \eIf{$\tilde \bv' > 0$}{Go to~$(**)$\; }{
\nl   set 
$\cS = {\rm supp}\, \tilde \bv'$, $X'
= \tilde X|_{\cS}$, and $\bv' = \tilde \bv|_{\cS}$. Go to~$(**)$.} }{
\nl we have $\cH \not \subset \cI$ and $\rho(\tilde X) <\rho(X)$.
Set  $X^{(k+1)} = \tilde X$ 
and go to the next $(k+1)$st iteration~(*)\; }}}}
\nl \Return{$\widehat X$ is a solution\;}
\caption{Algorithm~1 for minimising 
the spectral radius over an $L_1$-ball \label{algo1}}
\end{algorithm}

\bigskip

\begin{center}
\textbf{3.5. Explanations and proofs for Algorithm~1}
\end{center}
\bigskip

\noindent \textbf{Explanation for Algorithm~1}. 

The algorithm is a relaxation scheme: the value of the spectral radius 
$\rho(X)$ is non-increasing during the whole algorithm.
 
The algorithm consists of finitely many iterations (*), each iteration consists of 
several steps (**). In one iteration, during all steps except for the last one the value  $\rho(X)$
is the same. After the last step this value 
becomes strictly smaller, then the iteration is completed and we start the next iteration.  
\smallskip 

(*) In $k$th iteration we have a matrix 
$X = X^{(k)}$. We find its minimal leading eigenvector~$\bv$ and 
denote $\cS = {\rm supp}\, \bv$. Then till the end of this iteration we work 
 on the set $\cS$ without involving other indices. 
We consider the matrix $X'$, which is a restriction of $X$ to the set $\cS$. 
Then the vector $\bv' = \bv|_{\cS}$ is a positive minimal leading eigenvector 
of $X'$.  We find the basic set $\cH$ of $X'$ (Definition~\ref{d.20}), 
see~(\ref{eq.SH}): 

\begin{equation}\label{eq.SH}
X^{(k)} \ = \  
\begin{array}{l}
\overbrace{{\hspace{14mm}} }^{\cS}\hspace{13mm} \\ 
\underbrace{
\begin{array}{|ccc|ccc|}
\hline 
{}& {} & {} & {} & {} & {} \\
{}& X & {} & {} & * & {} \\
{}& {} & {} & {} & {} & {} \\
\hline 
{}& {} & {} & {} & {} & {} \\
{}& 0 & {} & {} & * & {} \\
{}& {} & {} & {} & {} & {} \\
\hline 
\end{array}}_{\Omega}
\\
{}
\end{array}
\ ; \qquad 
X \ = \ 
\begin{array}{c}
\overbrace{
\begin{array}{|ccc|ccc|}
\hline 
{}& {} & {} & {} & {} & {} \\
{}& B & {} & {} & * & {} \\
{}& {} & {} & {} & {} & {} \\
\hline 
{}& {} & {} & {} & {} & {} \\
{}& 0 & {} & {} & C & {} \\
{}& {} & {} & {} & {} & {} \\
\hline 
\end{array}}^{\cS}\\
\hspace{14mm} \underbrace{{\hspace{13mm}} }_{\cH}
\end{array}\, . 
\end{equation}

(**) The input of this loop is  a triple $(X, \bv, \cS)$, 
where $X$ is Boolean $n\times n$ matrix, $\bv$ is a minimal 
leading eigenvector of~$X$ and $\cS$ is the support of~$\bv$. The output 
is a new triple $(X, \bv, \cS)$, where $X$ is a modified matrix of the same size 
with the 
same set $\cS$ or  with a new (smaller) set $\cS$.

We denote $X' = X|_{\cS}, \, \bv' = \bv|_{\cS}$ and replace all rows of $X'$ by the minimal rows 
in the corresponding balls $\cB(A_i', r)$ with respect to~$\bv'$. 
All other elements of $X$ are not changed. 
We obtain a matrix $\tilde X$. 

If all rows of~$X'$ are already minimal, then $X'$ is a solution, Algorithm~1 terminates.  

If all rows with indices in~$\cH$, where $\cH \subset \cS$ is  
the basic set of~$X'$, are already minimal, then the leading eigenvalue of $\tilde X'$ is simple (Theorem~\ref{th.20} below).

We compute the minimal leading eigenvector~$\tilde \bv'$ of  $\tilde X'$. 
If $\tilde \bv'>0$, then  denote 
$\bv' = \tilde \bv', \, \tilde X' = X'$, keep the same $\cS$ and go to the next step (**). 
Otherwise we set $\cS = {\rm supp}\, \tilde \bv' = \cS$, restrict everything to 
this set and go to   (**). Thus, we have a sequence of matrices 
$X'$ with simple leading eigenvalues, unless the support $\cS$ gets smaller. 
Then we pass to the smaller support (i.e., to a submatrix), again obtain a sequence of matrices with simple leading eigenvelue, unless the support gets smaller, etc. 
We do it until in some step not all rows with indices in~$\cH$ are minimal. 
In this case  
$\rho(X') < \rho(X)$. We set $X^{(k+1)} = \tilde X$ and go to the next iteration (*).

\smallskip

In Algorithm~1 we have three main components:
\smallskip

\noindent  1) {\em Invariants.} In each iteration we have a Boolean matrix $X$,
its minimal leading eigenvector~$\bv$ and a set of indices
$\cS = {\rm supp }\, \bv$;
\smallskip

\noindent  2)  {\em Progress measure.} The spectral radius 
$\rho(X)$ strictly decreases in iterations.  Inside each iteration 
the spectral radius  is the same and 
the index set $\cS$ is non-increasing. When at some step $\rho(X)$ strictly decreases, we recompute the set $\cS$ and start the new iteration with this set. 
 
 Inside one step of the algorithm (in the inner loop) the algorithms does not cycle 
 (Theorem~\ref{th.20})
\smallskip

\noindent  3)  {\em Stopping criterion.} The algorithm stops when the current 
 matrix~$X' = X|_{\cS}$ is minimal in every row. In this case 
 $X$ is also minimal in every row and hence (Theorem~\ref{th.10})
 $X$ is a solution of problem~(\ref{eq.prob20}). 
\smallskip

\medskip

Now we are going to show that Algorithm 1 is well-defined and 
always finds a solution within finite time. 
The well-definedness means that at each iteration
matrix~$X'$ has a leading eigenvector, which is unique up
to normalization. We are proving more: $X'$ has a
simple leading eigenvalue. The claim that the algorithm 
finds the global solution in a finite time 
means two things: 
1) all statements formulated in the description of the  algorithm are correct;  2) the algorithm does not cycle. 

\begin{theorem}\label{th.20}
Algorithm 1 is well-defined. It finds the global solution
of problem~(\ref{eq.prob20}) in a finite number of steps.
\end{theorem}
{\tt Proof.} First we need to prove the correctness of interim conclusions: 
in steps 2 and~5 the matrix $X$ is optimal; 
in step 7 the assertion $\cH \subset \cJ$ implies that 
$\lambda_{\max}$ is simple and is not changed after this step, and 
in step 11, $\cH \not \subset \cJ$ implies that 
$\lambda_{\max}$ becomes stricly smaller.  Algorithm~1 is a modification of the algorithm from~\cite{NP2} 
derived specially for the uncertainty sets~$\cF_i = \cB(A_i, r)$. 
The proofs for the steps 3,5,7,11
are the same 
as the proof of analogous Theorem~8 from~\cite{NP2}. 
We only replace solutions of the problem 
$\la X_i, \bv \ra \to \min, \ X_i \in \cB(A_i, r)$ by 
the minimal rows of the matrix $X$ with respect to $\bv$ (Lemma~\ref{l.10}) and use the criterion of the solution of problem~(\ref{eq.prob20}) from Corollary~\ref{c.10}.
The proof of non-cyclicity is also the same as for Theorem~8  from~\cite{NP2}. 
We also note that due to Theorem~\ref{th.10}, Algorithm~1 runs over 
extreme points of the sets $\cB(A_i, r)$, i.e., over Boolean vectors. 
Since the total number of Boolean vectors is finite, the non-cyclicity 
implies that Algorithm~1 terminates within finite time. In step 2 the algorithm terminates since for an integer matrix $X$, $\rho(X)< 1$ means that $\rho(X) = 0$, 
hence $X$ has the minimal possible spectral radius.

{\hfill $\Box$}
\smallskip

\begin{center}
\textbf{4. Algorithmic solution of the max-MAS problem and generalizations}
\end{center}
\bigskip

\textbf{Solution of the max-MAS problem}. Take $r_0 \,  = \, \|A_i\|_{\infty} $. 
Since the ball $\cB(A, r_0)$ contains the zero matrix, it follows that 
$\min_{X \in \cB(A,r_0)}\rho(X) = 0$. Then applying Algorithm~1
and the integer bisection on the segment~$[0,r_0]$ we find the 
smallest $r$ such that $\min_{X \in \cB(A,r_0)}\rho(X) = 0$. 
For this~$r$, the matrix~$\widehat X$ provided by Algorithm~1 and 
the graph of this matrix 
give the answer. \textbf{Solution is  completed}. 
\medskip

Before looking at numerical results and discussing the complexity, 
we note that several generalisations of the max-MAS problem 
can be solved using slight modifications of our method. Below we formulate three of them.  
\medskip 

\noindent \textbf{Problem 1}. {\em To each vertex of a graph~$G$, a non-negative integer is assigned. Make the graph acyclic  by cutting 
at most the assigned number of incoming edges from each vertex.}
\medskip 

Let a number $r_i \ge 0$ be assigned to the $i$th vertex, $i = 1, \ldots , n$. 
Problem~1 is solved with Algorithm~1 by replacing all 
balls $\cB(A_i, r)$ with $\cB(A_i, r_i)$. 
If the minimal spectral radius is zero, then the answer is 
affirmative and the matrix of the desired graph is available. 

  \smallskip

\noindent \textbf{Problem 2}. {\em Solve the max-MAS problem for a weighted 
graphs, with given positive weights of edges.}
\medskip 

Solved by usual (non-integer) bisection and invoking 
Algorithm~1. One needs only to 
modify  the definition of minimal row as follows: 
the ordering of the $i$th row by numbers $v_j$ is replaced by ordering 
by numbers $\alpha_{ij}v_j$, where $\alpha_{ij}$ is the weight
of an edge from $g_j$ to $g_i$. 
  \smallskip

\noindent \textbf{Problem 3}. {\em Solve the max-MAS problem 
with an extra assumption that some of edges are ``untouchable'', 
i.e., it is prohibited to cut them. }
\medskip 

Solved as the usual max-MAX with the following 
modification of the definition of minimal row (Definition~\ref{d.20}):  
For each $i = 1, \ldots , n$, in $i$th row the positions of untouchable edges are removed from the set $\, {\rm supp}\, A_i$. 

\bigskip

\begin{center}
\textbf{5. The max-MAS problem and an approximate solution \\ for the classical MAS 
problem}
\end{center}
\bigskip

Having solved the max-MAS problem we obtain a matrix~$\widehat X$ and the 
corresponding acyclic graph, denote it by 
$G_0$. This graph can be considered as an approximate solution for the 
MAS problem for the graph~$G$. However,   usually $G_0$ has much less than $\frac12 |\cE|$ edges and so gives a bad approximation for MAS. The reason is obvious: 
the algorithm of solving the max-MAS problem tries to cut the maximal allowed number of incoming edges from each vertex 
and therefore cuts more edges than needed. Nevertheless, the following modified scheme gives satisfactory results: 
\smallskip 

\begin{algorithm}
\KwData{A graph $G$}
\KwResult{An acyclic subgraph $\bar G$, which is an approximate solution to 
MAS problem.}
\Begin{Apply Algorithm~1 to the graph~$G$. Obtain a solution 
$G_0$ to the max-MAS problem; } 

\nl Take the matrix $\widehat X$ of $G_0$. Find its Frobenius factorisation: 
$P^{-1}\widehat X P$, where $P$ is a permutation matrix \; 

\nl Set $Y_{ij} = \bigl[ P^{-1}AP\bigr]_{ij}$ if $j> i$ and $Y_{ij} = 0$ otherwise \;

\nl \Return{Set $\bar X = PYP^{-1}$. Then the graph $\bar G$ of the matrix $\bar X$ is an approximate solution for MAS\;}
\caption{Algorithm~2 for approximate solution of MAS \label{algo2}}
\end{algorithm}
\bigskip 

\noindent \textbf{Explanation for Algorithm~2}.  The solution $\widehat X$ of max-MAS problem has the spectral radius equal to zero. Hence its Frobenius factorisation $Z = P^{-1}\widehat X P$ 
is upper triangular with zero diagonal. Replacing the over-diagonal 
part of $Z$ by the over-diagonal part of the matrix $P^{-1}AP$ we keep the spectral radius equal to zero and reduce the distance to the matrix $P^{-1}AP$.
Denote the obtained matrix by $Y$. 
We have $\rho(Y)  = 0$ and the inverse permutation $PYP^{-1}$ is an approximate solution for MAS.  
\medskip 

\begin{remark}\label{r.30}{\em In fact Algorithm~2 finds the ordering of 
vertices in $\cV$ corresponding to the max-MAS solution~$G_0$. 
In this ordering (given by the permutation matrix~$P$) the matrix 
$\widehat X$ has an upper triangular form with zero diagonal. Then we set $\bar G$
to be the acyclic graph corresponding to this enumeration. 
}
\end{remark}

Note that Algorithm~2 can easily be modified to find approximate solutions of 
several generalizations of the MAS problem that inspired by Problems~1-3 
in Section~4. For example: 
\smallskip 

 {\em Find the maximal acyclic subgraph under the extra assumptions that at most $r_i$ incoming edges are cut from the $i$th vertex, 
$i=1, \ldots , n$, where $\{r_i\}_{i=1}^n$ are given integers.} 
\smallskip 

This corresponds to 
Problem~2. Problem~3 rises another variant of MAS: 
\smallskip 

{\em Find the maximal acyclic subgraph under the extra assumptions that some edges are untouchable.}

\smallskip

We are not aware of any known algorithms from the literature for approximate solutions of those problems.

\bigskip 

\begin{center}
\textbf{6. Numerical efficiency and complexity issue \\
for the algorithmic solution of the MAS problem}
\end{center}
\bigskip

In this section we demonstrate the practical efficiency 
of our methods. We show and discuss the results of numerical 
experiments with random graphs of various densities both for the MAX and the max-MAS 
problems. Then we discuss the complexity issue. 
\bigskip

\begin{center}
\textbf{6.1 Numerical results for random graphs of various density}
\end{center}
\bigskip

The algorithm  for the max-MAS problem demonstrates a very good efficiency. It consists in integer bisection in parameter 
$r$, where in each iteration of the bisection we solve problem~(\ref{eq.prob20})
with Algorithm~1. 
The total number of iterations therefore does not exceed $1+ \log_2  \|A\|_{\infty}$; in each step we apply Algorithm~1. 
Numerical results are shown in Tables~1 and~2. 
Table~1 demonstrates results for random graphs with sparsity  $\frac{|\cE|}{n^2}$ 
between $49\%$ and $91\%$, Table~2 shows results for sparsity 
between $5\%$ and $74\%$. For each
dimension  $n = |\cV|$ from $20$ to $1500$, we made $20$ experiments 
and put the average number of steps ($\#$ {\em steps}) 
and the average running time. Let us recall that by one step we mean 
one computation of the leading eigenvector, because this is the most expensive operation.  
The numerical experiments were performed on a standard laptop with the following specifications: Dell XPS 13 with Intel Core i7-6500U CPU @ 2.50GHz and 8GB RAM. The algorithm was coded in Python, and the code for the algorithm can be found on https://github.com/ringechpil/thesis_codes. 

Every time we use the obtained
solution for the max-MAS problem to find the approximate solution 
for the MAS problem (for the same graph). The rate of approximation 
$\gamma$ is written in the last row.  

In all our examples  Algorithm~1 finds
the solution within 3-5 steps and this number grows very slow with the dimension. 
Then to solve the max-MAS problems we need to apply Algorithm~1 at most 
$\log_2  n + 1$ times. We see form the Tables 1 and 2 that for graphs with 
$250$ vertices the complete solution of max-MAS problem 
takes  less than~$35$ steps, which is done for less than~$9$ seconds; for graphs with $1000$
vertices the solution takes less than~$11$ minutes. The average rate of approximation~$\gamma$ is quite stable and stays close to~$0.6$ for all dimensions. 
\smallskip

\begin{table}[H]
\begin{center}
\begin{tabular}{c|c c c c c}
$n$ & 50 & 250 & 500 & 1000 & 1500\\
\hline
$\mbox{time}$ & 0.36s & 8.1s & 66.42s & 622.43s & 2860.79s\\
$\# \, \mbox{steps}$ & 17 & 34.6 & 38.5 & 44.7 & 50.3\\
$\gamma$ & 0.644 & 0.621 & 0.615 & 0.616 & 0.616 

\end{tabular}
\caption*{{\footnotesize Table 1. Solving the max-MAS and approximating
 MAS for random graphs with sparsity $9-51\%$}}
\end{center}
\end{table}

\begin{table}[H]
\begin{center}
\begin{tabular}{c|c c c c c}
$n$ & 50  & 250 & 500 & 1000 & 1500\\
\hline
$\mbox{time}$ & 0.35s & 6.56s & 61.06s & 605.73s & 2614.02s\\
$\# \, \mbox{steps}$ & 18.9 & 32.1 & 41.8 & 43.1 & 43.7\\
$\gamma$ & 0.6 & 0.592 & 0.592 & 0.593 & 0.592 

\end{tabular}
\caption*{{\footnotesize Table 2. Solving the max-MAS and approximating
 MAS for random graphs with sparsity $26-95\%$}}
\end{center}
\end{table}   
    
   \bigskip

   \begin{center}
\textbf{6.2 The complexity of Algorithm~1 and comparison with other algorithms}
\end{center}
\bigskip
   
    The theoretical complexity of the max-MAS problem is not known for us. 
  We can only conjecture that  it is polynomial. It was shown in~\cite{CP} that 
  for  positive strictly convex smooth sets~$\cF_i$, the greedy algorithm has a quadratic convergence. This certainly explains the fast convergence of Algorithm~1 for finite sets $\cF_i$ but does not give good estimate for the theoretical complexity. 
  At least, an idea to approximate the balls $\cB(A_i, r)$ by convex smooth set
  does not work, since the parameters of quadratic convergence depend 
  on radii of curvature of the sets~$\cF_i$ which are too large for a 
  tight approximation. 

There are many algorithms in the literature  for approximate solving of MAS, 
see~\cite{CMM, HR} for short surveys and the bibliography.  
To the best of our knowledge, the most effective ones have the approximation factor 
$\delta \, = \, 0.5 \, + \, \Omega(\frac{\alpha}{\log n})$, 
where $0.5 + \alpha \, = \, \frac{|\widehat \cE|}{|\cE|}$ is the fraction of the original edges contained in the maximal acyclic subgraph~\cite{CMM}. 
Let us recall that for positive functions $f, g$, the symbol 
$f = \Omega(g)$ means that there is a constant $C > 0$ for which 
$\Omega (g(t)) \ge C\, g(t)$ for all~$t$.   The numerical results presented here show that our approach outputs graphs with $\delta \ge \gamma\approx 0.6$,  even for the large size $n$ or large density. Moreover, this estimate seems not to decrease with the growth of~$n$. 
Recall that $\gamma$ is a lower bound for the approximating factor 
$\delta$, where we replace $|\widehat \cE|$ by $|\cE|$. Therefore, 
the true values of $\delta$ in our examples are better (and can be much better)
than those given in the tables. Therefore, we see that our Algorithm~1 gives, at least 
in the numerical examples, not worse approximation than those known in the literature and 
performs very fast even for relatively large graphs.  

We also remark that the estimates given in Tables 1 and 2  are for randomly generated graphs with the only restriction of sparsity. As we 
we shall see in the next section, for graphs from applied problems
such as  small-world  networks  we get even better results.

\bigskip 

\begin{center}
\textbf{7. Applications}
\end{center}
\bigskip

Making the graph acyclic has a wide variety of applications. 
One of them  is discovering hierarchies within a graph \cite{Tat, Sun}. This is very useful for social networks, both real-world and virtual. Some other applications are for testing the electrical circuits \cite{Eve}, for ``telling stories'' \cite{Acu}, which is again used in biology for metabolic networks that describe biochemical road maps \cite{Tol}. Approximating the MAS can also be used for optimizing data flows/pipelines \cite{Kra} machine learning and artificial intelligence \cite{Koll, Liu}.

In \cite{Bri} a connection between a positive linear switching system and its asymptotic stability is established: a positive LSS is asymptotically stable if and only if its corresponding graph is acyclic. Using this fact and utilising our max-MAS algorithm along with an LSS stabilization algorithm introduced in \cite{Cve}, we can cut some interdependencies and construct a stable LSS from an unstable one, while keeping close to its original structure.  

We now apply our algorithm on small-world networks. 
Newman–-Watts-–Strogatz small-world graphs can be used for social network modelling \cite{Watt, Arn} (but also for modelling the networks in biology, epidemiology and neuroscience; see \cite{Hum} and references therein). They can be defined using three parameters: $n$ -- number of nodes forming a ring; $k$ -- the degree of each node, where each edge connects a node with its $k$ nearest neighbours; and $p$ -- the rewiring probability, i.e. a probability that an edge will be rewired from a neighbouring node to some random distant node. It is instructive to have a degree $k \gg \ln\, n$, but still not too large, in order not to make a graph overly dense and connected. Also, $p$ should not be too large, since rewiring too many edges makes the graph loose its small-world structure, more resembling a random network.

In the following series of numerical experiments we keep the parameters $k = 25, p = 0.1$ fixed, while we vary the number of vertices $n$. The results are shown in Table 3 below:\\

\begin{table}[H]
\begin{center}
\begin{tabular}{c|c c c c c}
$n$ & 50 & 250 & 500 & 1000\\
\hline
$\mbox{time}$ & 0.39s & 27.39s & 169.79s & 1607.81s\\
$\# \, \mbox{steps}$ & 13.6 & 43 & 53.8 & 73.4\\
$\gamma$ & 0.61 & 0.617 & 0.618 & 0.627

\end{tabular}
\caption*{{\footnotesize Table 3. Solving the max-MAS and approximating
 MAS for Newman--Watts--Strogatz small-world graphs with $k = 25,\ p = 0.1$.}}
\end{center}
\end{table}

Now we keep the number of vertices and rewiring $n = 500, p = 0.1$ fixed, while we vary the
 number of connected nearest neighbours $k$:

\begin{table}[H]
\begin{center}
\begin{tabular}{c|c c c c c}
$k$ & 5 & 10 & 25 & 100 & 250\\
\hline
$\mbox{time}$ & 76.07s & 76.98s & 126.63s & 138.15s & 119.65s\\
$\# \, \mbox{steps}$ & 69.8 & 39.2 & 55 & 37.6 & 30.6\\
$\gamma$ & 0.79 & 0.672 & 0.621 & 0.544 & 0.523 

\end{tabular}
\caption*{{\footnotesize Table 4. Solving the max-MAS and approximating
 MAS for Newman-Watts-Strogatz small-world graphs with $d = 500,\ p = 0.1$.}}
\end{center}
\end{table}

We see that as the network gets denser, our algorithm tends to cut significantly more edges. In this manner, our algorithm works better for sparser graphs, which is convenient, since small-world networks are usually not dense. We also perform tests for the graphs with a fixed number of vertices and node degree $n = 500, k = 25$, and we vary the rewiring probability $p$:

\begin{table}[H]
\begin{center}
\begin{tabular}{c|c c c c c}
$p$ & 0.02 & 0.3 & 0.6\\
\hline
$\mbox{time}$ & 167.51s & 80.54s & 73.94s\\
$\# \, \mbox{steps}$ & 71 & 25.6 & 22.6\\
$\gamma$ & 0.613 & 0.62 & 0.623 

\end{tabular}
\caption*{{\footnotesize Table 5. Solving the max-MAS and approximating
 MAS for Newman-Watts-Strogatz small-world graphs with $d = 500,\ k = 25$.}}
\end{center}
\end{table}

When it comes to the number of preserved edges, we notice no big changes as the graphs structure breaks more towards the random network.

\bigskip 

\noindent \textbf{Acknowledgements.} 
The authors are grateful to the 
anonymous Referees for attentive reading and for many valuable remarks.


\bigskip 

\begin{thebibliography}{NN}

\bibitem{Acu}
V.\,Acu\~na,  
 E.\,Birmel\'e, 
 L.\,Cottret, P.\,Crescenzi,
 F.\,Jourdanf, V.\,Lacroixa,
 A.\,Marchetti-Spaccamela, 
 A.\,Marinoe, 
P.V.\,Milreu, M.-F.\,Sagot, L.\,Stougie, 
\newblock {\em Telling stories: Enumerating maximal directed acyclic graphs with a constrained set of sources and targets}, 
\newblock  Theor. Comput. Sci.,  457 (2012), 1--9 . 
\smallskip

\bibitem{Akian1}
M.\,Akian, S.\,Gaubert, J.\,Grand-Clément, and J.\,Guillaud, 
\newblock {\em  The operator approach to entropy games}, 
\newblock  Theory Comput. Syst. 63 (2019), no 5, 1089-–1130
\smallskip 


\bibitem{An}
J.\,Anderson,
\newblock {\em Distance to the nearest stable Metzler matrix}
\newblock in: 2017 IEEE 56th Annual Conference
on Decision and Control (CDC), 2017, pp. 6567-–6572. 
\smallskip

\bibitem{Arn}
V.\,Arnaboldi, A.\,Passarella, M.\,Conti, and R.I.M.\,Dunbar,
\newblock {\em Online social networks: human cognitive constraints in facebook and twitter personal graphs}, 
\newblock  A volume in Computer Science Reviews and Trends, Elsevier (2015)

\bibitem{Asarin}  
E.\,Asarin, J.\,Cervelle, A.\,Degorre, C.\,Dima, F.\,Horn, and V.\,Kozyakin, 
\newblock {\em   Entropy games
and matrix multiplication games}, 
\newblock  In 33rd Symposium on Theoretical Aspects of Computer
Science (STACS 2016), February 2016, Orl\'eans, France, 11:1 -– 11:14. 
\smallskip 


\bibitem{BP}
A.\,Berman and R.J.\,Plemmons,
\newblock {\em Nonnegative matrices in the mathematical sciences},
\newblock  Academic Press, Now York, 1979.
\smallskip

\bibitem{BS}
B.\,Berger, P.W.\,Shor, 
\newblock {\em Tight bounds for the Maximum Acyclic Subgraph
problem},
\newblock  J. Algorithms, 25 (1997), no 1, 1--18.  
\smallskip 

\bibitem{Bri}
C.\,Briat,
\newblock {\em Sign properties of Metzler matrices with applications},
\newblock Linear Alg. Appl.,  515 (2017), 53--86.  
\smallskip


\bibitem{By}
R.\,Byers,
\newblock {\em A bisection method for measuring the distance of a stable to unstable matrices},
\newblock SIAM J. on
Scientific and Statistical Computing, 9 (1988),  875--881.
\smallskip

\bibitem{CMM}
M.\,Charikar, K.\,Makarychev, and Yu.\,Makarychev, 
\newblock {\em On the advantage over
random for maximum acyclic subgraph},
\newblock  48th Annual IEEE Symposium
on Foundations of Computer Science (FOCS'07) (2007)
\smallskip 


\bibitem{Cve}
A.\,Cvetkovi\'c,
\newblock {\em Sign properties of Metzler matrices with applications},
\newblock Calcolo 57, (2020), no 1. 
\smallskip

\bibitem{CP}
A.\,Cvetkovi\'c and V.Yu.\,Protasov,
\newblock {\em The greedy strategy for optimising the Perron eigenvalue},
\newblock  submitted to Mathematical Programming, 	arXiv:1807.05099 (2018). 
\smallskip

\bibitem{Eve}
G.\,Even, J.\,Naor, B.\,Schieber, M.\,Sudan, 
\newblock {\em Approximating minimum feedback sets and multicuts in directed graphs}, 
\newblock  Journal of Electronic Testing,  1 (1990), 163-–174.  
\smallskip 

\bibitem{G}
F.R.\,Gantmacher, 
\newblock {\em The theory of matrices},  
\newblock Chelsea, New York, 2013.
\smallskip 


\bibitem{GS}
N.\,Gillis and P.\,Sharma, 
\newblock {\em On computing the distance to stability for matrices using linear
dissipative Hamiltonian systems},
\newblock Automatica 85 (2017), 113--121.
\smallskip


\bibitem{GM}
N.\,Guglielmi and M.\,Manetta.
\newblock {\em Approximating real stability radii},
\newblock IMA Journal of Numerical Analysis, 35 (2015), no. 3, 1402--1425.
\smallskip


\bibitem{GP2}
N.\,Guglielmi and V.\,Protasov, 
\newblock {\em On the closest stable/unstable nonnegative matrix
and related stability radii},  
\newblock SIAM J. Matrix Anal. 39 (2018), no 4,  1642--1669. 
\smallskip 

\bibitem{GMR}
V.\,Guruswami, R.\,Manokaran, and P.\,Raghavendr, 
\newblock {\em Beating the random ordering is hard: inapproximability of maximum acyclic subgraph}, 
\newblock FOCS 2009, Conference proceedings. 
\smallskip  



\bibitem{HR}
R.\,Hassin and S.\,Rubinstein, 
\newblock {\em Approximations for the maximum acyclic subgraph problem}, 
Information Processing Letters,  51 (1994) no 3, 133--140. 


\bibitem{HK}  
A.J.\,Hoffman and R.M.\,Karp, 
\newblock {\em  On nonterminating stochastic games}, 
\newblock  Management Science.
Journal of the Institute of Management Science. Application and Theory Series, 12 (1966), 359-–370.
\smallskip 

\bibitem{Hum}
M. D. Humphries and K Gurney,
\newblock {\em Network ``Small-World-Ness'': A quantitative method for determining canonical network equivalence}, 
\newblock PLOS ONE, 3 (2008), no 4.  

\bibitem{K}
R.M.\,Karp, 
\newblock {\em Reducibility among combinatorial problems. In Complexity of Computer
Computations}, 
\newblock New York: Plenum (1972),  85–-103.
\smallskip 


\bibitem{Koll}
D.\,Koller and N.\,Friedman,
\newblock {\em Probabilistic Graphical Models: Principles and Techniques}, 
\newblock  Adaptive Computation and Machine Learning series, The MIT Press (2009)
\smallskip 

\bibitem{Kra}
R.\,Kramer, R.\,Gupta, and M.L.\,Soffa,
\newblock {\em The combining DAG: a technique for parallel data flow analysis}, 
\newblock IEEE Trans. Paral. Distr. Systems 5 (1994), no 8, 805--813. 
\smallskip 

\bibitem{Liu}
Z.\,Liu, K.\,Li and X.\,He, 
\newblock {\em Cutting cycles of conditional preference networks with feedback set approach},
\newblock Computational Intelligence and Neuroscience (2018)
\smallskip 

\bibitem{MMS}
C.~Mehl, V.~Mehrmann, and P.~Sharma,
\newblock {\em Stability radii for linear
hamiltonian systems with dissipation under structure-preserving
perturbations,}
\newblock  SIAM Journal on Matr. Anal. Appl.,
 37 (2016), no. 4,  1625--1654.
\smallskip

\bibitem{NP1}
Y.Nesterov and V.Yu.Protasov,
\newblock {\em Optimizing the spectral radius},
\newblock SIAM J. Matrix Anal. Appl. 34 (2013), no 3,
999--1013.
\smallskip

\bibitem{NP2}
Y.\,Nesterov and V.Yu.\,Protasov,
\newblock {\em Computing closest stable
non-negative matrix}, 
\newblock SIAM J. Matrix Anal. Appl.,   41 (2020), no 1, 1–-28.
\smallskip

\bibitem{New}
A.\,Newman,  
\newblock {\em The maximum acyclic subgraph problem and degree-3 graphs}, 
\newblock Approximation, randomization, and combinatorial optimization (Berkeley, CA, 2001), 147-–158, Lecture Notes in Comput. Sci., 2129, Springer, Berlin, 2001.
\smallskip 

\bibitem{ONVD}
F.X.~Orbandexivry, Y.~Nesterov, and P.~Van Dooren,
\newblock {\em Nearest stable system using successive
convex approximations},
\newblock  Automatica, 49 (2013), 1195--1203.
\smallskip

\bibitem{PY}
I.Post and Y.Ye, 
\newblock {\em The simplex method is strongly polynomial for deterministic markov decision processes}, 
\newblock Mathematics of Operations Research,  published online (2015), 
https://doi.org/10.1287/moor.2014.0699

\bibitem{P16}
V.Yu.\,Protasov,
\newblock {\em The spectral simplex method},
\newblock Math. Prog., 156 (2016), 485--511
\smallskip

\bibitem{R}  
U.G.\,Rothblum, 
\newblock {\em Multiplicative Markov decision chains}, 
 \newblock Math. Oper. Research,
9 (1984), 6-–24.


\bibitem{Sun}
J.\,Sun, D.\,Ajwani, P.K.\,Nicholson, A.\,Sala, S.\,Parthasarathy, 
\newblock {\em Breaking Cycles In Noisy Hierarchies}, 
\newblock  WebSci '17: Proceedings of the 2017 ACM on Web Science Conference, 151–-160 (2017). 
\smallskip 

\bibitem{Tat}
N.\,Tatti,
\newblock {\em Faster Way to Agony}, 
\newblock  In: Calders T., Esposito F., H\"ullermeier E., Meo R. (eds) Machine Learning and Knowledge Discovery in Databases. ECML PKDD 2014. Lecture Notes in Computer Science, vol 8726. Springer, Berlin, Heidelberg (2014)
\smallskip

\bibitem{Tol}
J.B.\,Toledo et al,
\newblock {\em Metabolic network failures in Alzheimer's disease: A biochemical road map}, 
\newblock  Alzheimers Dement. 9 (2017), 965--984.  
\smallskip 


\bibitem{T}
R.\,Tarjan, 
\newblock {\em Depth first search and linear graph algorithms}, 
\newblock SIAM J. Comput. 1 (2) (1972),  146-–160.
\smallskip 


\bibitem{Watt}
D.\,Watts and  S.\,Strogatz,
\newblock {\em Collective dynamics of ``small-world'' networks}, 
\newblock  Nature, 393 (1998), 440-–442. 

\end{thebibliography}
\end{document}